\documentclass[12pt,draft]{article}
\usepackage{amssymb,amsmath,amsfonts,curves,ascmac,fancybox,makeidx,color,
amsthm,amsxtra,bm,cases}

\allowdisplaybreaks

\newcommand{\tr}{{\rm{tr}}}

\newtheorem{thm}{Theorem}       
\newtheorem{prop}{Proposition}  
\newtheorem{lem}{Lemma}         
\newtheorem{cor}{Corollary}     

\title{Linear relations among roots of a polynomial}
\author{Yoshiyuki Kitaoka}
\date{}
\begin{document}
\maketitle

Let $f(x)\in\mathbb Z[x]$ be a monic polynomial of degree $n$ with   roots $\alpha=\alpha_1,\dots,\alpha_n\in\mathbb C$, and let us consider a linear relation among roots
\begin{equation}
\label{eq1}
m_1\alpha_1+\dots+m_n\alpha_n=m\quad(m_1,\dots,m_n,m\in\mathbb Q).
\end{equation}
By $\alpha_1+\dots+\alpha_n= \tr(f)\in \mathbb Z$, there is always such a linear relation.
If $m_1=\dots=m_n$ holds at \eqref{eq1}, then we call it a trivial linear relation, otherwise non-trivial.
A trivial linear relation  is reduced to
$\alpha_1+\dots+\alpha_n= \tr(f)$, and
we pointed out the importance of the linear relations among roots when we study the distribution  of local roots of $f(x)\equiv 0$ modulo a prime (\cite{K1},\cite{K2}). 

It is known in \cite{K1} that if $f(x)$ is irreducible and  if the degree $n$ is prime or the Galois group of $f(x)$ is isomorphic to $S_n,A_n$,
then there is no non-trivial  linear relation among roots,
and that in the case of $n=4$ there is a non-trivial linear relation if and only if $f(x)=g(h(x))$ for some quadratic polynomials $g,h$, and in the case that $n=6$ and the Galois group is abelian, the polynomials with non-trivial linear relations among roots are classified. 
But even the case that the Galois group is isomorphic to $S_3$ was incomplete in \cite{K1}.

In this paper, we study the problem group-theoretically, and 
throughout this paper, we assume that $f(x)$ is a monic irreducible polynomial over the rational number field $\mathbb{Q}$ with a complex root $\alpha$, and  
$\mathbb{Q}(\alpha)/\mathbb{Q}$ is a Galois extension with Galois group
 $G:=Gal(\mathbb{Q}(\alpha)/\mathbb{Q})$.

We consider the group algebra $\mathbb{C}[G]$ and generalize the action of $G$ on the root $\alpha$
 to $\mathbb{C}[G]$ as $M[\alpha]:=\sum_{g\in G} m_gg(\alpha)
\in\mathbb{C}$ 
for $M=\sum_{g\in G}m_gg\in \mathbb{C}[G]$.
In case of $M\in G$, $M[\alpha]$ is the ordinary action of the element $M\in G$ on $\alpha$,
and if $M\in \mathbb{Q}[G]$, then $M[\alpha]\in\mathbb{Q}(\alpha)$.
We define the vector spaces $LR,LR_0$ of linear relations  by
\begin{align}\label{eq2}
LR &:=\{(M,m)\in \mathbb{Q}[G]\times\mathbb{Q}\mid m=M[\alpha] \in \mathbb{Q}\},
\\ \label{eq3}
LR_0 &:=\{M\in\mathbb{Q}[G]\mid M[\alpha]\in\mathbb{Q}\}.
\end{align}

It is known (\cite{L}) that  by making $G$ act on $\mathbb{C}[G]$ from the left,
$\mathbb{C}[G]$ contains the trivial representation with multiplicity $1$ and every  irreducible representation of $G$, that is
by putting 
$$
\delta:=|G|^{-1}\sum_{g\in G}g,
$$
the $G$-stable subspace $V$ spanned by $g-\delta$ $(g\in G)$ contains all non-trivial irreducible representations   and does not contain the trivial representation.
 The equations $\delta^2=\delta$, $g\delta=\delta g=\delta$ for $g\in G$, $1-\delta=-\sum_{g\not=1}(g-\delta)$, and $\delta V=V\delta=0$
 are clear where $1$ is the unit of the group $G$.
It is also obvious that $\mathbb{C}[G]=\mathbb{C}\delta\oplus V$, and
the decomposition $\sum_g m_gg=(\sum_g m_g)\delta+\sum_g m_g(g-\delta)$ implies
\begin{align*}
LR_0=\mathbb{Q}\delta\oplus (V\cap LR_0).
\end{align*}
We say that an element $M\in \mathbb{C}[G]$ is trivial if and only if $M\in\mathbb{C}\delta$.
And
there is a non-trivial linear relation among roots if and only if $\dim LR_0>1$.
%
\begin{prop}\label{prop1}
For   $M\in\mathbb{C}[G]$,   $M$ is trivial 
 if and only if   $M$ annihilates $V$.
 \end{prop}
\proof
{
The ``only if''-part follows from $\delta V=0$.
Conversely, we suppose that  $MV=0$ with $M=\sum m_gg$, in particular $M(1-\delta)=0$ for the unit $1\in G$; then
\begin{align*}
M(1-\delta)&=\sum_g m_{g}(g - \delta)
\\
&=\sum_{g\ne1} (m_{g}-m_{1})(g - \delta).
\end{align*}
Thus we have $m_{g} = m_{1}$ for every $g\in G$, since the set of $g-\delta$ $(g\ne1)$ is a basis of $V$.
Therefore $M\in\mathbb{C}\delta$ follows.
\qed
}

Let us  note the following fundamental equation for $M\in\mathbb{Q}[G]$
\begin{align}\label{eq4}
\sum_{g\in G}g(M[\alpha])g^{-1}=M(\sum_{g\in G}g(\alpha)g^{-1}),
\end{align}
because, putting $M=\sum_{h\in G}m_hh$ with $m_h \in \mathbb{Q}$,
\begin{align}\nonumber
\sum_{g\in G}g(M[\alpha])g^{-1}&=
\sum_{g\in G}(\sum_{h\in G}m_hgh(\alpha))(ghh^{-1})^{-1}
\\\nonumber
&=(\sum_{h\in G}m_hh)(\sum_{g\in G}g(\alpha)g^{-1}).
\end{align}
%
\begin{prop}\label{prop2}
Let $M\in\mathbb{Q}[G]$; then
 $M\in LR_0$ holds  if and only if
 \begin{equation}\label{eq5}
M(\sum_{g\in G}g(\alpha)g^{-1})V=0.
 \end{equation}
\end{prop}
\proof{
The \eqref{eq5} is equivalent to $\sum_{g\in G}g(M[\alpha])g^{-1}V=0$ by \eqref{eq4}.
Suppose that  $M\in LR_0$; then $m:=M[\alpha]\in\mathbb{Q}$ and 
 it is easy to see $\sum_{g\in G}g(M[\alpha])g^{-1}\linebreak[3]
=\sum_{g\in G}g(m)g^{-1}=m|G|\delta$,
and it annihilates $V$.
Conversely, suppose that \eqref{eq5}; then  Proposition \ref{prop1} implies that
 $\sum_{g\in G} g(M[\alpha])g^{-1}$ is trivial, that is $M[\alpha]\in\mathbb{Q}(\alpha)$ is
 fixed by $G=Gal(\mathbb{Q}(\alpha)/\mathbb{Q})$, i.e., $M[\alpha]$ is a rational number, 
 hence $M\in LR_0$.
\qed
}

\vspace{1mm}
Let $\chi$ be an irreducible character of $G$ and 
denote the central idempotent in $\mathbb{C}[G]$ to $\chi$ by
$$
c_\chi:=\chi(1)|G|^{-1} \sum_{g\in G}\overline{\chi(g)}g.
$$
They satisfy
$$
\sum_\chi c_\chi=1,\,\,c_{\chi_1}c_{\chi_2}=\delta_{\chi_1,\chi_2}c_{\chi_1}.
$$
If $\chi$ is the trivial character, then $c_\chi=\delta$. 
And $V$ is the direct sum of $\mathbb{C}[G]c_\chi$ with non-trivial irreducible characters $\chi$.
If $A$ is a matrix representation corresponding to $\chi$, then $\mathbb{C}[G]c_\chi$ is a direct sum of minimal ideals, which give the representation $A$.

If $\chi$ is an irreducible character, then the mapping $g\mapsto \sigma(\chi(g))$ is also an irreducible character for every $ \sigma\in Gal(\mathbb{C}/\mathbb{Q})$.
We denote it by $\sigma(\chi)$,
and let $\Psi$ be a minimal  set of irreducible characters such that every irreducible character is of the form
$\sigma(\chi)$ for ${}^\exists \sigma\in Gal(\mathbb{C}/\mathbb{Q}),{}^\exists \chi\in \Psi$. 

For an irreducible character $\chi$, we put
 $$
\mathbb{Q}(\chi):=\mathbb{Q}(\{\chi(g)\mid g \in G\}),
$$
which is an abelian extension over $\mathbb{Q}$
and we  denote the Galois group by $Gal(\chi):=Gal(\mathbb{Q}(\chi)/\mathbb{Q})$. 
Let $K$ be  a splitting field  of  $G$ such that $K$ is  a Galois extension of $\mathbb{Q}$,
and $\mathbb{Q}(\chi)\subset K$.
For example, $\mathbb{Q}(\zeta_e)$ with exponent $e$ of $G$ is such a splitting field, denoting a primitive $n$th root of unity by $\zeta_n$.

We put
$$
C(\chi) := \sum_{\sigma \in Gal(\chi)} c_{\sigma(\chi)}
=\chi(1)|G|^{-1}\sum_{g\in G} \tr_{\mathbb{Q}(\chi)/\mathbb{Q}}(\chi(g))g\,\,(\,\in\mathbb{Q}[G]\,),
$$
which is a central idempotent satisfying $\sum_{\chi\in\Psi}C(\chi)=1$.
We note that $C(1)=\delta$ and
 \eqref{eq5} is equivalent to
$$
M(\sum_{g\in G}g(\alpha)g^{-1})C(\chi)=0
$$
for every $\chi(\ne1)\in\Psi$.

For an irreducible character $\chi$, we put
\begin{equation}\label{eq6}
LR(\chi):=\{M\in  \mathbb{Q}[G]   \mid M=MC(\chi),
 M(\sum_gg(\alpha)g^{-1})=0
\}.
 \end{equation}
 If $\chi $ is the trivial character, then
$$
LR(1)=\{m\delta\mid m\in\mathbb{Q}, m\,\tr_{\mathbb{Q}(\alpha)/\mathbb{Q}}(\alpha)=0\}.
$$
 We note that 
 for $M\in\mathbb{Q}[G]$ with  $MC(\chi)=M$,
the condition $ M(\sum_gg(\alpha)g^{-1})\linebreak[3]=0$ is equivalent to 
$$ M(\sum_gg(\alpha)g^{-1})C(\chi)\mathbb{Q}[G]=0,$$
i.e., by defining a matrix representation  $A$ by $g(v_1,\dots,v_t)=(v_1,\dots,v_t)A(g)$
for a basis $\{v_1,\dots,v_t\}$ of $C(\chi)\mathbb{Q}[G]$ over $\mathbb{Q}$, 
the condition above is equivalent to
$\sum_{h,g\in G}m_hg^{-1}(\alpha)A(hg)=0$, putting $M=\sum_{h\in G}m_hh$.
The rational representation $A$ is decomposed to the sum of representations corresponding to
characters $\sigma(\chi)$.
 %
 %
%
 \begin{prop}\label{prop3}
For an irreducible non-trivial character $\chi$, we have
\begin{equation}\label{eq7}
LR(\chi)=\{MC(\chi)\mid M\in LR_0\}\,(=LR_0C(\chi)),
\end{equation}
and we have a linear relation
\begin{equation}\label{eq8}
M[\alpha]=0\,\,\text{ for }M\in LR(\chi).
\end{equation}
 \end{prop}
\proof
{
Let $M\in LR(\chi)$; then $M=MC(\chi)$ is clear by definition.
The fundamental equation  \eqref{eq4} implies  $\sum_gg(M[\alpha])g^{-1} =0$, 
that is $g(M[\alpha]) =0$ for every $g\in G$,
hence we have
$M[\alpha]=0$, i.e., \eqref{eq8} and hence $M \in LR_0$.
Conversely, let $M\in LR_0$; then $MC(\chi)\in\mathbb{Q}[G]$
and $(MC(\chi))C(\chi) = MC(\chi)$ are  obvious.
Proposition \ref{prop2} implies $M(\sum g(\alpha)g^{-1})V\linebreak[3]=0$,
which implies $MC(\chi)(\sum g(\alpha)g^{-1})=0$ by $C(\chi)\delta=0$, i.e., $MC(\chi) \in LR(\chi)$.
\qed
}
%
\begin{cor}\label{cor1}
We have
\begin{align}\label{eq9}
V\cap  LR_0 =\oplus_{\chi(\ne1)\in\Psi} LR(\chi),
\end{align}
\end{cor}
\proof
{
Let $M\in V\cap  LR_0$;
then  we see that $M\delta=0$  implies $M=M(1-\delta)=\sum_{\chi(\ne1)\in\Psi}MC(\chi)$.
Proposition \ref{prop3} implies $MC(\chi)\in LR(\chi)$.
Hence the left-hand side is contained in the right-hand side of \eqref{eq9}.
Conversely, suppose that $M\in LR(\chi)$ $(\chi\ne1)$; then we see that  $M=MC(\chi)\in V$ and   \eqref{eq8} implies  $M\in LR_0$.
Thus the right-hand side of \eqref{eq9} is contained in the left-hand side. 
\qed
}
%
\begin{cor}\label{cor2}
A linear relation among roots is a sum of a trivial linear relation and $M_\chi(\alpha)=0$ for several
 $M_\chi\in LR(\chi)$ with  $\chi(\ne1)\in\Psi$.
\end{cor}
\proof
{
Let $(M,m)\in LR$. 
Decompose as $M=\sum m_gg=(\sum m_g)\delta+M'$ with $M':=\sum m_g(g-\delta)$. 
Then $M'\in V \cap  LR_0$ is clear.
Corollary \ref{cor1} and \eqref{eq8} complete the proof.
\qed
}

\vspace{1mm}
From now on, we take and  fix an irreducible character $\chi$,
a minimal ideal $\mathfrak{m}$ of a simple algebra $K[G]c_\chi $ with
a basis  $\{w_1,\dots,w_{\chi(1)}\}$ of $\mathfrak{m}$ over $K$, and define an irreducible matrix representation $A$ corresponding to $\chi$ by
$$
g(w_1,\dots,w_{\chi(1)})=(w_1,\dots,w_{\chi(1)})A(g).
$$
The correspondence $g\mapsto A(g)$ is extended to  an isomorphism from 
$\mathbb{C}[G]c_\chi$ to the matrix ring $M_{\chi(1)}(\mathbb{C})$  with 
$A(K[G]c_\chi)=M_{\chi(1)}(K)$ because of $K$ being a splitting field.
Through this correspondence, take $v_{i,j}\in K[G]c_\chi $ satisfying
 \begin{equation*}
v_{i,j}(w_1,\dots,w_{\chi(1)})=(w_1,\dots,w_{\chi(1)})E_{j,i}, \text{ i.e., }v_{i,j}w_k=\delta_{i,k}w_j,
\end{equation*}
where $E_{j,i}$ is a matrix whose $(j,i) $-entry is $1$, otherwise $0$.
Thus, the set $\{v_{i,j}\mid 1\le i,j\le \chi(1)\} $ is a basis of $K[G]c_\chi$ over $K$.

If $\deg\chi=1$, i.e., $\chi$ is a homomorphism from $G$ to $\mathbb{C}$, then $K[G]c_\chi=Kc_\chi$ and we can take $c_\chi$ as  $w_1,v_{1,1}$.

In the following, the property of the field $K$ which we use is that $K/\mathbb{Q}$ is a Galois extension with $\mathbb{Q}(\chi)\subset K$ and $K[G]c_\chi$ is isomorphic to $M_{\chi(1)}(K)$. 
%
%
\begin{lem}\label{lem1}
We have
\begin{gather}\label{eq10}
g(v_{i,1},\dots,v_{i,\chi(1)})=(v_{i,1},\dots,v_{i,\chi(1)})A(g)\text{ for } {}^\forall g\in G,
\\\label{eq11}
v_{i,j}v_{k,l}=\delta_{i,l}v_{k,j}\text{ for } 1\le i,j,k,l\le\chi(1). 
\end{gather}
\end{lem}
\proof
{
Write $gv_{i,j}=\sum_{1\le k,l\le\chi(1)}x_{k,l}v_{k,l}$ with $x_{k,l}\in K$;
then we see that $gv_{i,j}w_m$ is equal to $\sum_{k,l}x_{k,l}v_{k,l}w_m=
\sum_{k,l}x_{k,l}\delta_{k,m}w_l=\sum_{1\le l\le\chi(1)}x_{m,l}w_l  =\sum_{1\le k\le\chi(1)}x_{m,k}w_k $.
On the other hand,  we  find that $gv_{i,j}w_m\linebreak[3]=g\delta_{i,m}w_j=\delta_{i,m}\sum_{1\le k\le\chi(1)}w_kA(g)_{k,j}$,
where $A(g)_{k,j}$ is the $(k,j)$-entry of $A(g)$.
Comparing the coefficients of $w_k$ of two equations, we have $x_{m,k}=\delta_{i,m}A(g)_{k,j}$.
Hence  we have $gv_{i,j}=\sum_{1\le m,k\le\chi(1)}\delta_{i,m}A(g)_{k,j}v_{m,k}= \sum_{1\le k\le\chi(1)}A(g)_{k,j}v_{i,k} $, i.e., \eqref{eq10}.
The second equation follows from
\begin{align*}
v_{i,j}v_{k,l}(w_1,\dots,w_{\chi(1)})&=v_{i,j}(w_1,\dots,w_{\chi(1)})E_{l,k}
\\
&=(w_1,\dots,w_{\chi(1)})E_{j,i}E_{l,k}
\\
&=\delta_{i,l}(w_1,\dots,w_{\chi(1)})E_{j,k}
\\
&=\delta_{i,l}v_{k,j}(w_1,\dots,w_{\chi(1)}).
\end{align*}

\qed
}

Let us recall that the mapping ${M}\mapsto A({M})$ from $\mathbb{C}[G]c_\chi$ to 
$M_{\chi(1)}(\mathbb{C})$ defined by
$$
{M}(w_1,\dots,w_{\chi(1)})=(w_1,\dots,w_{\chi(1)})A({M})
$$
is an isomorphism.
\begin{lem}\label{lem2}
For $M=\sum_{1\le k,l\le\chi(1)}m_{k,l}v_{l,k}\in K[G]c_\chi$ with $m_{k,l}\in K$, 
we have   
\begin{equation}\label{eq12}
A(M(\sum_{g\in G}g^{-1}(\alpha)g))=(m_{i,j})(\sum_{g\in G}g^{-1}(\alpha)A(g)).
\end{equation}
\end{lem}
\proof
{

It is easy to see that
\begin{align*}
M(w_1,\dots,w_{\chi(1)})&=\sum_{k,l}m_{k,l}v_{l,k}(w_1,\dots,w_{\chi(1)})
\\
&=\sum_{k,l}m_{k,l}(w_1,\dots,w_{\chi(1)})E_{k,l}
\\
&=
(w_1,\dots,w_{\chi(1)})(m_{i,j}),
\end{align*}
and
$$
\sum_{g\in G}g^{-1}(\alpha)g(w_1,\dots,w_{\chi(1)})=
(w_1,\dots,w_{\chi(1)})(\sum_gg^{-1}(\alpha)A(g)),
$$
hence
$$
M(\sum_{g\in G}g^{-1}(\alpha)g)(w_1,\dots,w_{\chi(1)})=
(w_1,\dots,w_{\chi(1)})(m_{i,j})(\sum_gg^{-1}(\alpha)A(g)),
$$
i.e., \eqref{eq12}.
\qed
}
\begin{prop}\label{prop4}
Suppose that $K\cap\mathbb{Q}(\alpha)=\mathbb{Q}$; then $\sum_{g\in G}g^{-1}(\alpha)A(g)$
is the zero matrix if and only if  the first column of $\sum_{g\in G}g^{-1}(\alpha)A(g)$ is  the zero vector.  
\end{prop}
\proof
{
Suppose that  $K\cap\mathbb{Q}(\alpha)=\mathbb{Q}$;
then taking elements $c_{j,g}\in K$ such that $w_j = \sum_{h\in G} c_{j,h}hw_1$,
we have $gw_j=\sum_{h\in G}c_{j,h}ghw_1$, hence
\begin{align*}
gw_j&=(w_1,\dots,w_{\chi(1)})\times(\text{\,the $j$th column of } A(g))
\\
&=(w_1,\dots,w_{\chi(1)})\times(\sum_h c_{j,h}\times(\text{\,the first column of } A(gh))),
\end{align*}
therefore
\begin{align*}
\text{\,the $j$th column of } A(g)
=\sum_h c_{j,h}\times(\text{\,the first column of } A(gh)).
\end{align*}
Thus we have 
\begin{align*}
&\sum_{g\in G} g^{-1}(\alpha)\times(\text{ the $j$th column of }A(g))
\\
=\,& \sum_{g,h\in G} g^{-1}(\alpha)c_{j,h}\times(\text{ the first column of }A(gh)) 
\\
=\,& \sum_{g,h\in G} c_{j,h}hg^{-1}(\alpha)\times(\text{ the first column of }A(g)) 
\\
=\,& \sum_{h\in G} c_{j,h}h\left(\sum_{g\in G}g^{-1}(\alpha)\times(\text{ the first column of }A(g)) \right),
\end{align*}
where $h\in G=Gal(\mathbb{Q}(\alpha)/\mathbb{Q})$ is extended to $Gal(\mathbb{C}/\mathbb{Q})$ with
$h$ being the identity on $K$.
If, therefore the first column of $\sum_{g\in G}g^{-1}(\alpha)A(g)$ is  the zero vector,
then  $\sum_{g\in G}g^{-1}(\alpha)A(g)$ is the zero matrix.
The converse is obvious.
\qed
}

\noindent
{\bf Remark}  The condition $\sum_{g\in G}g^{-1}(\alpha)A(g)=0$ can occur at least in the case
of $\deg \chi=1$ as in the next section.
What is a necessary and/or sufficient condition to it?

When we make $\sigma\in Gal(\mathbb{C}/\mathbb{Q})$ act on $\mathbb{C}[G]$ by
$\sum_gm_gg\mapsto \sum_g\sigma(m_g)g$, we use $\sigma^\circ$, i.e.,
$\sigma^\circ(\sum_gm_gg ):=\sum_g\sigma(m_g)g$.
Then we see that $\sigma^\circ(\mathfrak{m})$ is a minimal ideal of $K[G]c_{\sigma(\chi) }$ and 
the set $\{\sigma^\circ(w_1),\dots,\linebreak[3]\sigma^\circ(w_{\chi(1)})\}$ is
a basis  of the ideal $\sigma^\circ(\mathfrak{m})$ over $K$, and
$$
g(\sigma^\circ(w_1),\dots,\sigma^\circ(w_{\chi(1)}))
=(\sigma^\circ(w_1),\dots,\sigma^\circ(w_{\chi(1)}))\sigma(A(g)).
$$
and
$$
\sigma^\circ(v_{i,j})(\sigma^\circ(w_1),\dots,\sigma^\circ(w_{\chi(1)}))\linebreak[3]=(\sigma^\circ(w_1),\dots,\sigma^\circ(w_{\chi(1)}))E_{j,i}.
$$

%
%
%
\begin{prop}\label{prop5}
We define a mapping $\psi$ from $\{M\in \mathbb{Q}[G]\mid M=MC(\chi)\}$ to
$\mathbb{Q}(\chi)[G]c_\chi$  by $\psi(M)= Mc_\chi$; then $\psi$  is isomorphic.
For $M=MC(\chi)\in\mathbb{Q}[G]$, write 
$$
Mc_\chi=\sum_{1\le k,l\le\chi(1)}m_{k,l}v_{l,k}\in K[G]c_\chi\quad (m_{k,l}\in K).
$$
Then $M\in LR(\chi)$ holds  if and only if
\begin{equation}\label{eq13}
\sigma((m_{i,j}))\left(\sum_{g\in G} g^{-1}(\alpha)\sigma(A(g))\right)=0^{(\chi(1))}\,\text{ for }\,{}^\forall\sigma\in Gal(K/\mathbb{Q}).
\end{equation}
If, in particular $\mathbb{Q}(\alpha)\cap K=\mathbb{Q}$, then   $M\in LR(\chi)$
is equivalent to
\begin{equation}\label{eq14}
(m_{i,j})\times(\text{ the first column of }\sum_{g\in G} g^{-1}(\alpha)A(g))
=0^{(\chi(1),1)}.
\end{equation}
\end{prop}
\proof{
Since $C(\chi)$ is the sum of $c_{\sigma(\chi)}$ with  $\sigma\in Gal(\mathbb{Q}(\chi)/\mathbb{Q})$,
the  injectivity of $\psi$ follows from  $M=MC(\chi)=\sum_{\sigma\in Gal(\chi)}Mc_{\sigma(\chi)}$.
Conversely for a given element $\hat N=N_0c_\chi \in\mathbb{Q}(\chi)[G]c_\chi$ $(N_0\in \mathbb{Q}(\chi)[G])$,
 $ N:=\sum_{\sigma\in Gal(\chi)}\sigma(N_0)c_{\sigma(\chi)}\in\mathbb{Q}[G]$ satisfies  $ N=NC(\chi)$ because of  
$c_{\sigma(\chi)}C(\chi)= c_{\sigma(\chi)}   $.
By $Nc_\chi=\hat N$,  the mapping  $\psi$ is an isomorphism.
We find that  $MC(\chi)(\sum_gg^{-1}(\alpha)g)=0$ if and only if $MC(\chi)(\sum_gg^{-1}(\alpha)g)c_{\sigma(\chi)}=0$ for  $M\in\mathbb{C}[G]$ and $\sigma\in Gal(\mathbb{Q}(\chi)/\mathbb{Q})$.
Hence for $M=MC(\chi)\in \mathbb{Q}[G]$, the condition $M\in LR(\chi)$ is equivalent to   
$M(\sum_gg^{-1}(\alpha)g)c_{\sigma(\chi)}=\sigma^\circ(Mc_{\chi})(\sum_gg^{-1}(\alpha)g)=0$ for every $\sigma\in Gal(K/\mathbb{Q})$.
Applying Lemma \ref{lem2} to $\sigma^\circ(Mc_\chi)=\sum\sigma(m_{k,l})\sigma^\circ(v_{l,k})\in K[G]c_{\sigma(\chi)}$ instead of $M$ there,
we find that the equation is $\sigma((m_{i,j}))(\sum_g g^{-1}(\alpha)\sigma(A(g)))=0^{(\chi(1))}$,
i.e., \eqref{eq13}.
Suppose that $\mathbb{Q}(\alpha)\cap K=\mathbb{Q}$; then every 
$\sigma\in Gal(K/\mathbb{Q})$ is extended to $Gal(\mathbb{C}/\mathbb{Q})$ with $\sigma(\alpha)=\alpha$,
hence 
\eqref{eq13} is equivalent to
\begin{equation}\label{eq15}
(m_{i,j})\sum_{g\in G} g^{-1}(\alpha)A(g)=0^{(\chi(1))}.
\end{equation}
As in the proof of Proposition \ref{prop4}, take   $c_{j,g}\in K$ such that $w_j = \sum_{h\in G} c_{j,h}hw_1$.
Then we have 
\begin{align*}
&\sum_{g\in G} g^{-1}(\alpha)\times(\text{ the $j$-column of }A(g))
\\
=\,& \sum_{g,h\in G} c_{j,h}hg^{-1}(\alpha)\times(\text{ the first column of }A(g)) 
\\
=\,& \sum_{h\in G} c_{j,h}h\left(\sum_{g\in G}g^{-1}(\alpha)\times(\text{ the first column of }A(g)) \right),
\end{align*}
where we extend  the automorphism $h\in G=Gal(\mathbb{Q}(\alpha)/\mathbb{Q})$ to $\mathbb{C}$
so that $h$ is the identity on $K$.
Now \eqref{eq14} follows from \eqref{eq15}. 
\qed
}
\vspace{1mm}

\noindent
Suppose that  we can take all $v_{l,k}$ in $\mathbb{Q}(\chi)[G]c_\chi$;
then all coefficients $m_{k,l}$ defined by 
 $Mc_\chi=\sum_{1\le k,l\le\chi(1)}m_{k,l}v_{l,k}$
$ (m_{k,l}\in K)$
 are in $\mathbb{Q}(\chi)$,
hence the mapping  $\psi$ induces an isomorphism from 
$\{M\in \mathbb{Q}[G]\mid M=MC(\chi)\}$ to
$M_{\chi(1)}(\mathbb{Q}(\chi))$.

\vspace{1mm}
\noindent
{\bf Remark} The proof shows that if $\mathbb{Q}(\chi)= \mathbb{Q}$ 
and $v_{i,j}\in\mathbb{Q}[G]$ for ${}^\forall i,j$,
then $M\in LR(\chi)$ if and only if \eqref{eq14} holds.

Does the equivalence of conditions $M\in LR(\chi)$ and  \eqref{eq14}  hold unconditionally?
%
%
%
%
%
%
%
\section{ The case of degree $\bold 1$} 
In this section, we study the case of $\deg\chi=1$.
\begin{thm}\label{th1}
Suppose that the degree of an irreducible character $\chi$ is $1$ and let $l$ be the order  of $\chi$;
then we have $LR(\chi)=\{0\}$ or $\sum_{g\in G} \chi(g)^dg(\alpha)=0 $ for every integer $d$ relatively prime to $l$.
If $LR(\chi)\ne\{0\}$, i.e., $\sum_{g\in G} \chi(g)^dg(\alpha)=0 $ holds for every integer $d$ relatively prime to $l$,
then we have
\begin{align*}
LR(\chi)
 &=\{\sum_g \tr_{\mathbb{Q}(\zeta_l)/\mathbb{Q}}(\,m\,\overline{\chi(g)}\,) g \mid m \in \mathbb{Q}(\zeta_l)\}
\end{align*}
with $\dim_{\mathbb{Q}} LR(\chi)=[\mathbb{Q}(\zeta_l):\mathbb{Q}]$.
\end{thm}
\proof
{
For each element $h\in G$, the equation $hc_\chi=\chi(h)c_\chi$ is clear.
We note $\mathbb{Q}(\chi)=\mathbb{Q}(\zeta_l)$, and take $c_\chi$ as $w_1,v_{1,1}$ as noted.
Let $M=MC(\chi)\in\mathbb{Q}[G]$ and $M_0:=Mc_\chi=mc_\chi$ with $m\in \mathbb{Q}(\chi)$; then $M=\sum_{\sigma\in Gal(\chi)}Mc_{\sigma(\chi)}=
\sum_{\sigma\in Gal(\chi)}\sigma^\circ(M_0)=|G|^{-1}\sum_gtr_{\mathbb{Q}(\chi)/\mathbb{Q}}
(m\overline{\chi(g)})g$ is clear.
We parametrize $M$ by $m$ as in Proposition \ref{prop5}.
By \eqref{eq13}, $M\in LR(\chi)$ holds if and only if $\sigma(m)\sum_g g^{-1}(\alpha)\sigma(\chi(g))=0$ for every $\sigma\in Gal(\chi)$,
which means either $m=0$ or $\sum_g g^{-1}(\alpha)\sigma(\chi(g))=0$.
Hence,  $LR(\chi)\ne0$ implies $\sum_g g^{-1}(\alpha)\sigma(\chi(g))=0$ for every $\sigma\in Gal(\chi)$, and then $M\in LR(\chi)$ holds for every $m\in\mathbb{Q}(\chi)$. 
Since $\sum_g \tr_{\mathbb{Q}(\zeta_l)/\mathbb{Q}}(\,m\,\overline{\chi(g)}\,) g=0$ if and only if $m =0$, we find $\dim_{\mathbb{Q}} LR(\chi)=[\mathbb{Q}(\zeta_l):\mathbb{Q}]$.
\qed
}

\vspace{5mm}

\noindent
{\bf Remark} 
It is not difficult to see that there are elements $a_i\in\mathbb{Q}(\zeta_l)$ such that $\chi(g)^d=
\sum_i \tr_{\mathbb{Q}(\zeta_l)/\mathbb{Q}}(a_i\overline{\chi(g)})\zeta_l^{i-1}$ $({}^\forall g\in G)$
 for  an integer $d$ relatively prime to  $l$, that is $LR(\chi)$ gets back the linear relation  $\sum_g \chi(g)^dg(\alpha)=0$ because of
  $\sum_g \chi(g)^dg(\alpha)=\sum_i\zeta_l^{i-1}(\sum_g\tr_{\mathbb{Q}(\zeta_l)/\mathbb{Q}}(a_i\overline{\chi(g)})g(\alpha))$.
%
%

\vspace{2mm}
We know 
$$
\tr_{\mathbb{Q}(\zeta_l)/\mathbb{Q}}(\zeta_l^k)=
[\varphi(l):\varphi(l/d)]\mu(l/d)\text{ for } d:=(k,l).
$$

\vspace{1mm}

\noindent
{\bf Remark} 
Suppose that the order $l$ of $\chi$ is prime; 
then by taking an element $g_0\in G$ with $\chi(g_0)=\zeta_l$ and by denoting the field 
fixed by $\ker \chi$ by $M_\chi$, we see that  $\sum \chi(g)^dg(\alpha)$ with $(d,l)=1$ is equal to
\begin{align*}
\sum_{k=0}^{l-1}\zeta_l^{dk} g_0^k(\tr_{\mathbb{Q}(\alpha)/M_\chi}(\alpha))&=\sum_{k=0}^{l-2}\zeta_l^{dk} (g_0^k( \tr_{\mathbb{Q}(\alpha)/M_\chi}(\alpha))
-g_0^{l-1}( \tr_{\mathbb{Q}(\alpha)/M_\chi}(\alpha)))
\end{align*}
by $1+\zeta^d+\dots+\zeta^{d(l-2)}+\zeta^{d(l-1)}=0$.
 If, hence $\sum\chi(g)^dg(\alpha)=0$ holds for every integer $d$ satisfying $( d,l)=1$,
then $g_0^k( \tr_{\mathbb{Q}(\alpha)/M_\chi}(\alpha))$ is independent of $k$,
therefore a linear relation among roots
 $ \tr_{\mathbb{Q}(\alpha)/M_\chi}(\alpha)= \tr(f)/l\in\mathbb Q$ follows.
 Here we used that the matrix $(\zeta_l^{dk})$ $((d,l)=1,0\le k\le l-2)$ is regular.

\vspace{2mm}
If $G=Gal(\mathbb{Q}(\alpha)/\mathbb{Q})$ is abelian, then the degree of every irreducible character
is $1$, hence Corollary \ref{cor2} and Theorem \ref{th1} describe linear relations among roots.

\noindent
{\bf Example 1}
Let $G=\langle g_0\rangle$ be of prime order $l$.
If $LR(\chi)\ne0$ for some non-trivial irreducible character $\chi$, 
then $\ker \chi=\{1\}$, i.e., $M_\chi=\mathbb{Q}(\alpha)$ is clear and  we have
the contradiction $\alpha= \tr_{\mathbb{Q}(\alpha)/M_\chi}(\alpha)= \tr(f)/l\in\mathbb Q$
as Remark above.
Hence $LR(\chi)=0$ for every non-trivial irreducible character, that is  there is no non-trivial linear relation. 

This is true  in general for an irreducible polynomial of prime degree  (Proposition 2 in \cite{K1}).

\noindent
{\bf Example 2}
Let $G$ be of order $4$ and suppose $LR(\chi)\ne 0$ for a non-trivial irreducible character.
In case that the order of $\chi$ is $4$,  $G$ is cyclic and
we  number roots by $\alpha_i = g_0^i(\alpha)$ for a generator $g_0$ of $G$.
Then we see that $\chi(g_0)=\sqrt{-1}$ implies $\sum_g\chi(g)^dg(\alpha)=
\sqrt{-1}^d\alpha_1 +(-1)^d\alpha_2+(-\sqrt{-1})^{d}\alpha_3+\alpha_4=0$ for every odd integer $d$,
which implies a contradiction $\alpha_2=\alpha_4$.
Next, suppose that the order of $\chi$ is $2$; then as Remark above we have $\sum_{g\in\ker\chi}g(\alpha)=\tr(f)/2$.
Hence $f(x)$ is a polynomial in $x^2 - (  \tr(f)/2  )x$.
Thus, if there is a non-trivial linear relation among roots for an irreducible abelian polynomial of degree $4$, then a polynomial is decomposable.

This is true if $f(x)$ is irreducible and of degree $4$ (\cite{K1}).

\noindent
{\bf Example 3}
Let $G=\langle g_0\rangle$ be a cyclic group of order $6$ and let $\omega$ be a third root of unity, i.e.,  $\omega^2+\omega+1=0$.
We number roots by $\alpha_i = g_0^i(\alpha)$.
Let $\chi$ be a non-trivial character of order $l$ satisfying $LR(\chi)\ne 0$.
Then the relations corresponding to $LR(\chi)$ are $\sum_g \tr_{\mathbb{Q}(\chi(g_0))/\mathbb{Q}}(m\overline{\chi(g)})g(\alpha)=0$ for every $m\in\mathbb{Q}(\chi(g_0))$.

\noindent
Case of $l=6$: We assume $\chi(g_0)=-\omega$.
\!\!Equations $\sum_g \tr_{\mathbb{Q}(\omega)/\mathbb{Q}}(m\overline{\chi(g)})g(\alpha)\linebreak[3]=0$ for $m=1,\omega$ imply $\alpha_1+\alpha_5+2\alpha_6=\alpha_2+2\alpha_3+\alpha_4$ and
$2\alpha_1+\alpha_2+\alpha_6=\alpha_3+2\alpha_4+\alpha_5$,
hence $\alpha_1-\alpha_4=-\alpha_2+\alpha_5=\alpha_3-\alpha_6$.
Therefore we have $f(x) = \{(x-\alpha_1)(x-\alpha_4)\} \{(x-\alpha_2)(x-\alpha_5)\} \{(x-\alpha_3)(x-\alpha_6)\}$ and three quadratic factors are polynomials over the cubic subfield
with rational discriminant by $(\alpha_1-\alpha_4)^2=(\alpha_2-\alpha_5)^2=(\alpha_3-\alpha_6)^2$.
It is of type (c) in Proposition 6 in \cite{K1}.

\noindent
Case of $l=3$: We may assume $\chi(g_0)=\omega$.
Then relations are $\alpha_1+\alpha_2+\alpha_4+\alpha_5=2\alpha_3+2\alpha_6,
\alpha_2+\alpha_3+\alpha_5+\alpha_6=2\alpha_1+2\alpha_4$,
i.e., $\alpha_1+\alpha_4=\alpha_2+\alpha_5=\alpha_3+\alpha_6$.
we have $f(x) = \{(x-\alpha_1)(x-\alpha_4)\} \{(x-\alpha_2)(x-\alpha_5)\} \{(x-\alpha_3)(x-\alpha_6)\}$ and three factors are polynomials in $x^2-(\alpha_1+\alpha_4)x=x^2-(tr(f)/3)x\in \mathbb{Q}[x]$.
It is of type (d) or (f).

\noindent
Case of $l=2$: We have $\chi(g_0)=-1$.
Relations are $\alpha_1+\alpha_3+\alpha_5=\alpha_2+\alpha_4+\alpha_6$,
and 
 we have $f(x) = \{(x-\alpha_1)(x-\alpha_3)(x-\alpha_5)\}\{(x-\alpha_2) (x-\alpha_4)(x-\alpha_6)\}$.
It is of type (b) or (e).

Type (f) is a combination of  two characters $\chi_1,\chi_2$ with $\chi_1(g_0)=\omega$ and $\chi_2(g_0)=-1$,
because $\chi_1(g_0)=\omega$ implies $\alpha_6= \tr(f)/3 - \alpha_3,\alpha_5= \tr(f)/3 - \alpha_2, \alpha_4=\tr(f)/3 - \alpha_1$ and then $\chi_2(g_0)=-1$ implies $\alpha_1+\alpha_3+(\tr(f)/3-\alpha_2)=
\alpha_2+(\tr(f)/3-\alpha_1)+(\tr(f)/3-\alpha_3)$, hence $\alpha_1-\alpha_2+\alpha_3=\tr(f)/6$.
\vspace{1mm}

\noindent
\section{ The case of $\bold S_3$  }
In this section, we take an irreducible polynomial $f(x)$ of degree $6$ with root $\alpha$
such that $\mathbb{Q}(\alpha)$ is a Galois extension of $\mathbb{Q}$ with Galois group 
isomorphic to $S_3$.
Denote the Galois group $Gal(\mathbb{Q}(\alpha)/\mathbb{Q})$ by $G$, and
let $G=\langle\sigma,\mu\rangle$ with $\sigma^3=\mu^2=1,\mu\sigma\mu=\sigma^2$.
The representations of $G$ are
\begin{enumerate}
\item\label{c'1}
the trivial character $\chi_1$, and $c_{\chi_1}=\frac{1}{6}\sum_gg$,
\item\label{c'2}
the character $\chi_2$ of degree $1$ defined by $\chi(\sigma)=1,\chi(\mu)=-1$,
and $c_{\chi_2}=\sum_{i=0}^2\sigma^i -\sum_{i=0}^2\sigma^i\mu$,
\item\label{c'3}
the character of degree $2$ corresponding to the representation 
$$
A(\sigma)=\left(\begin{array}{cc}0&-1\\1&-1\end{array}\right),\,
A(\mu)=\left(\begin{array}{cc}0&1\\1&0\end{array}\right)
$$
with $c_{\chi_3}=\frac{1}{3}(2-\sigma-\sigma^2) \in \mathbb{Q}[G]$
and $\mathbb{Q}[G]c_{\chi_3}=\langle 1-\sigma,\sigma-\sigma^2,\mu-\sigma\mu,\sigma\mu-\sigma^2\mu\rangle$.
\end{enumerate}
In case of $\chi_2$, we have only to invoke Theorem \ref{th1}.
Let us consider the case $\chi:=\chi_3$.
For simplicity, we number roots of $f(x)$ as follows:
\begin{equation*}
\begin{array}{lll}
\alpha_1 := \alpha,&\alpha_2:=\sigma(\alpha),&\alpha_3:=\sigma^2(\alpha),
\\
\alpha_4:=\mu(\alpha),&\alpha_5:=\sigma(\alpha_4),&\alpha_6:=\sigma^2(\alpha_4),
\\
\mu(\alpha_2)=\alpha_6,&\mu(\alpha_3)=\alpha_5,&
\end{array}
\end{equation*}
and abbreviate as

\begin{equation}\nonumber
m_1=m_{\sigma^0},\,
m_2=m_{\sigma},\,m_3=m_{\sigma^2},\,m_4=m_{\mu},\,m_5=m_{\sigma\mu}
,\,m_6=m_{\sigma^2\mu}.
\end{equation}
Put
\begin{align*}
v_1&:=(1-\sigma-\sigma\mu+\sigma^2\mu)/3\,(=v_1c_\chi),
\\
v_2&:=-(\sigma-\sigma^2-\mu+\sigma\mu)/3\,(=v_2c_\chi=-v_1\sigma),
\end{align*}
then it is easy to see that
 for $i=1,2$
 \begin{align*}
\sigma(v_i,\sigma v_i)=(v_i,\sigma v_i)A(\sigma),\mu(v_i,\sigma v_i)=(v_i,\sigma v_i)A(\mu),
\end{align*}
and
\begin{align*}
v_1(v_i,\sigma v_i)=(v_i,\sigma v_i)\left(\begin{array}{rr}1&0\\0&0\end{array}\right),
\sigma v_1(v_i,\sigma v_i)=(v_i,\sigma v_i)\left(\begin{array}{rr}0&0\\1&0\end{array}\right),
\\
v_2(v_i,\sigma v_i)=(v_i,\sigma v_i)\left(\begin{array}{rr}0&1\\0&0\end{array}\right),
\sigma v_2(v_i,\sigma v_i)=(v_i,\sigma v_i)\left(\begin{array}{rr}0&0\\0&1\end{array}\right).
\end{align*}
Thus $v_{1,1}:=v_1,v_{1,2}:=\sigma v_1,v_{2,1}:=v_2,v_{2,2}:=\sigma v_2$ satisfy \eqref{eq10}
and \eqref{eq11}, and $v_{i,j}$ are in $\mathbb{Q}[G]$ for  $i,j=1,2$.
These  imply, for $M=c_1v_1+c_2\sigma v_1+c_3v_2+c_4\sigma v_2$
\begin{equation}\label{eq16}
M(v_i,\sigma v_i)=(v_i,\sigma v_i)\left(
\begin{array}{rr}c_1&c_3\\c_2&c_4\end{array}
\right)\,(=(v_i,\sigma v_i)C, \text{ say}).
\end{equation}
Moreover we see, for $i=1,2$
\begin{align*}
&(\sum_g g^{-1}(\alpha)g)(v_i,\sigma v_i)
\\=\,&(v_i,\sigma v_i)
\left(\begin{array}{cc}
\alpha_1-\alpha_2-\alpha_5+\alpha_6&\alpha_2-\alpha_3+\alpha_4-\alpha_6
\\
-\alpha_2+\alpha_3+\alpha_4-\alpha_5&\alpha_1-\alpha_3+\alpha_5-\alpha_6
\end{array}\right)
\\
=&\,(v_i,\sigma v_i){\bf{A}}, \text{ say},
\end{align*}
hence the condition $M\in LR(\chi)$ is equivalent to $C{\bf{A}}=0.$
The second column of ${\bf{A}}$ is the image of the first column by $\sigma$,
hence $C{\bf{A}}=0$ is equivalent to $C\times(\text{the first column of }{\bf{A}})=0$.
Let us see  $({\bf{A}}_{1,1},{\bf{A}}_{2,1})\ne(0,0)$.
If ${\bf{A}}_{1,1}={\bf{A}}_{2,1}=0$ holds,
then we get ${\bf{A}}_{1,1}+\sigma({\bf{A}}_{2,1})=2\alpha_1-\alpha_2-\alpha_3=0$ 
and $2\alpha_2-\alpha_3-\alpha_1=0$, acting $\sigma$. Their  difference is $3(\alpha_1-\alpha_2)=0$,  which is a contradiction.

\noindent
Let us see $LR(\chi)$ explicitly.
\begin{thm}\label{th2}
Suppose that  $LR(\chi)\ne\{0\}$;
then there are rational numbers $a,b$ with $(a,b)\ne(0,0)$ such that
\begin{equation}\label{eq17}
a(\alpha_1-\alpha_2-\alpha_5+\alpha_6)=
b(\alpha_2-\alpha_3-\alpha_4+\alpha_5),
\end{equation}
and the basis of $LR(\chi)$ is given by $av_1+bv_2,\sigma(av_1+bv_2)$,
and corresponding linear relations are spanned by
$a(\alpha_1-\alpha_2-\alpha_5+\alpha_6)+b(-\alpha_2+\alpha_3+\alpha_4-\alpha_5)=0$
and $a(\alpha_2-\alpha_3+\alpha_4-\alpha_6)+b(\alpha_1-\alpha_3+\alpha_5-\alpha_6)=0$.
\end{thm}
\proof
{
Suppose $LR(\chi)\ne0$ and take a non-zero element $M=
c_1v_1+c_2\sigma v_1+c_3v_2+c_4\sigma v_2\in \mathbb{Q}[G]C(\chi)$.
 Then at least one of $c_i$ is not zero,
 hence there are rational numbers $a,b$ with $(a,b)\ne(0,0)$ satisfying $a{\bf{A}}_{1,1}+b{\bf{A}}_{2,1}=0$,
 i.e., \eqref{eq17}.
 The we have $\left|\begin{array}{rr}{\bf{A}}_{1,1}&b\\-{\bf{A}}_{2,1}&a\end{array}\right|
 =\left|\begin{array}{rr}{\bf{A}}_{1,1}&c_3\\-{\bf{A}}_{2,1}&c_1\end{array}\right|=
 \left|\begin{array}{rr}{\bf{A}}_{1,1}&c_4\\-{\bf{A}}_{2,1}&c_2\end{array}\right|=0$
by $C{\bf{A}}=0$.
Therefore,   there are rational numbers $\kappa_1,\kappa_2$
 such that $(c_3,c_1)=\kappa_1(b,a),(c_4,c_2)=\kappa_2(b,a)$, thus we have
 $M= \kappa_1(av_1+bv_2)+\kappa_2(a\sigma v_1+b\sigma v_2)$.

 \qed
}
\vspace{2mm}

\noindent
{\bf Remark}
Put $\hat M=x(\alpha_1-\alpha_2-\alpha_5+\alpha_6)+y(-\alpha_2+\alpha_3+\alpha_4-\alpha_5)
+z(\alpha_2-\alpha_3+\alpha_4-\alpha_6)+w(\alpha_1-\alpha_3+\alpha_5-\alpha_6)$.
When $(x,y,z,w) = (2/3,2/3,1/3,1/3) $, the relation $\hat M=0$ gives $\alpha_1+\alpha_4=\sigma(\alpha_1+\alpha_4)$. 
Similarly, $(x,y,z,w)=(-1/3,2/3,-2/3,4/3), (2/3,-1/3,-2/3,1/3)$ give the relation $\alpha_1+\alpha_5=\sigma(\alpha_1+\alpha_5), \alpha_1+\alpha_6=\sigma(\alpha_1
+\alpha_6)$, respectively. 
These three cases mean that the trace of $\alpha$ to a cubic subfield is rational.
When $(x,y,z,w)=(0,0,-1,1) ,(1,0,0,0), (0,1,0,1)$
give the relation
$\alpha_1-\alpha_4=\sigma(\alpha_1-\alpha_4),\alpha_1-\alpha_5=\sigma(\alpha_1-\alpha_5),\alpha_1-\alpha_6=\sigma(\alpha_1-\alpha_6)$, respectively.
These three cases mean that the discriminant of the minimal quadratic polynomial of $\alpha$ over a cubic subfield is rational.  
These are equations pointed out in \cite{K1}.

\noindent
{\bf Example}
Let $F(x)=x^6+3$ and $\zeta$ be a root of $F(x)=0$.
For a primitive third root $\omega$ of unity, define the automorphisms $\sigma$, $\mu$ by $\sigma(\zeta)=\omega\zeta$ and $\mu(\zeta)=-\zeta$.
In the following table, the left is a root of polynomial of $f$ and the middle is $\alpha_1+\alpha_2+\alpha_3$ if it is rational, and the right is ${\bf A}_{2,1}/{\bf A}_{1,1}$ if it is a non-zero 
rational number, and $\infty$ is put if ${\bf A}_{1,1}=0$.
$$
\begin{array}{|c|c|c|}
\hline
\text{root of }f & \alpha_1+\alpha_2+\alpha_3 & {\bf A}_{2,1}/{\bf A}_{1,1}
\\\hline
\zeta&0&-1
\\\hline
-\zeta^2-\zeta^3+\zeta^5&&\infty
\\\hline
-\zeta-\zeta^3+\zeta^4&&2
\\\hline
\zeta^2+\zeta^3+\zeta^5&&0
\\\hline
\zeta+\zeta^3+\zeta^4&&1/2
\\\hline
-2\zeta^2-2\zeta^3+\zeta^5&&3
\\\hline
\zeta+\zeta^3+\zeta^5&&-1
\\\hline
\zeta^2+\zeta^3+\zeta^4&&1
\\\hline
5-5\zeta-3\zeta^4&15&2/7
\\
\hline
\end{array}
$$
\section{The case of the dihedral group of order $\bf8$}
This case is similar to the previous case.
Let  $G:=Gal(\mathbb{Q}(\alpha)/\mathbb{Q})$ be the dihedral group defined by the fundamental relations
$\sigma^4=\mu^2=1, \mu^{-1}\sigma\mu=\sigma^3$.
The irreducible characters and the corresponding central idempotents are
\begin{enumerate}
\item\label{c1}
the trivial character $\chi_1$, and $c_{\chi_1}=\frac{1}{8}\sum_g g$,
\item\label{c2}
the character $\chi_2$  of degree $1$ defined by $\chi_2(\sigma)=1,\chi_2(\mu)=-1$,
and $$c_{\chi_2}=\frac{1}{8}(\sum_{i=0}^3 \sigma^i -\sum_{i=0}\sigma^i\mu),$$
\item\label{c3}
the character $\chi_3$  of degree $1$ defined by $\chi_3(\sigma)=-1,\chi_3(\mu)=1$, and
$$
c_{\chi_3}=\frac{1}{8}(1-\sigma+\sigma^2-\sigma^3+\mu-\sigma\mu+\sigma^2\mu-\sigma^3\mu),
$$
\item\label{c4}
the character $\chi_4$  of degree $1$ defined by $\chi_4(\sigma)=-1,\chi_4(\mu)=-1$,
and
$$
c_{\chi_4}=\frac{1}{8}(1-\sigma+\sigma^2-\sigma^3-\mu+\sigma\mu-\sigma^2\mu+\sigma^3\mu),
$$
\item\label{c5}
the character $\chi_5$ corresponding to the representation defined by
$$
\sigma\mapsto \left(\begin{array}{cc}0&-1\\1&0\end{array}\right),\,\,
\mu\mapsto\left(\begin{array}{cc}0&1\\1&0\end{array}\right)
$$ and 
$$
c_{\chi_5}=\frac{1}{2}(1-\sigma^2).
$$
\end{enumerate}
All central idempotents $c_\chi$ are in $\mathbb{Q}[G]$,
and the cases of non-trivial characters of degree $1$ are done by Theorem \ref{th1} and the linear relations are
$c_{\chi_i}(\alpha) = 0$ for $i=2,3,4$ if exists.
%
Hereafter we are concerned with the fifth character $\chi_5$.
For simplicity, we put $\chi:=\chi_5$;
then we have
\begin{align*}
\mathbb{Q}[G]c_\chi&=\langle 1-\sigma^2,\sigma-\sigma^3,(1-\sigma^2)\mu,(\sigma-\sigma^3)\mu\rangle.
\end{align*}
Put
\begin{align*}
{v_1}&:=(1-\sigma^2-\sigma\mu+\sigma^3\mu)/4=(1-\sigma^2)(1-\sigma\mu)/4=v_1c_\chi,
\\
{v_2}&:=-(\sigma-\sigma^3-\mu+\sigma^2\mu)/4=-(1-\sigma^2)(\sigma-\mu)/4=-v_2c_\chi,
\end{align*}
then for $i=1,2$
\begin{equation*}
\sigma({v_i},\sigma {v_i})=({v_i},\sigma {v_i})\left(\begin{array}{rr}0&-1\\1&0\end{array}\right),
\,\,\mu({v_i},\sigma {v_i})=({v_i},\sigma {v_i})\left(\begin{array}{rr}0&1\\1&0\end{array}\right).
\end{equation*}
Moreover, we find 
\begin{align*}
{v_1}({v_i},\sigma {v_i})=({v_i},\sigma {v_i})\left(\begin{array}{rr}1&0\\0&0\end{array}\right)&,
\sigma {v_1}({v_i},\sigma {v_i})=({v_i},\sigma {v_i})\left(\begin{array}{rr}0&0\\1&0\end{array}\right),
\\
{v_2}({v_i},\sigma {v_i})=({v_i},\sigma {v_i})\left(\begin{array}{rr}0&1\\0&0\end{array}\right)&,
\sigma {v_2}({v_i},\sigma {v_i})=({v_i},\sigma {v_i})\left(\begin{array}{rr}0&0\\0&1\end{array}\right).
\end{align*}
Hence for $M=c_1{v_1}+c_2\sigma {v_1}+c_3{v_2}+c_4\sigma {v_2}\in LR(\chi)$, we find
$$
M({v_i},\sigma {v_i})=({v_i},\sigma {v_i})\left(\begin{array}{rr}c_1&c_3\\c_2&
c_4\end{array}\right)\quad(i=1,2).
$$
Furthermore,  we see that $\sum_g g^{-1}(\alpha)g({v_i},\sigma {v_i})=({v_i},\sigma {v_i}){\bf A}$ 
$(i=1,2)$ for
\begin{equation*}
{\bf A}=
\left(
\begin{array}{rr}
\alpha_1-\alpha_3-\alpha_6+\alpha_8&\alpha_2-\alpha_4+\alpha_5-\alpha_7
\\
-\alpha_2+\alpha_4+\alpha_5-\alpha_7&\alpha_1-\alpha_3+\alpha_6-\alpha_8
\end{array}
\right),
\end{equation*}
where we abbreviate as $\alpha_1=\alpha,\alpha_i=\sigma^{i-1}(\alpha_1),(i=2,3,4)$,
and $\alpha_i=\sigma^{i-5}\mu(\alpha_1)$ $(i=5,6,7,8)$.
Hence the condition $M\in LR(\chi)$ is equivalent to
\begin{equation*}
\left(\begin{array}{rr}c_1&c_3\\c_2&c_4\end{array}\right)
\left(\begin{array}{r}{\bf A}_{1,1}\\{\bf A}_{2,1}\end{array}\right)=\left(\begin{array}{r}0\\0\end{array}\right),
\end{equation*}
since the second column of the matrix $\bf A$ is the image of the first column by $\sigma$.
The equation above is equivalent to
$$
\left|\begin{array}{rr}c_1&-{\bf A}_{2,1}\\c_3&{\bf A}_{1,1}\end{array}\right|=
\left|\begin{array}{rr}c_2&-{\bf A}_{2,1}\\c_4&{\bf A}_{1,1}\end{array}\right|=0.
$$
If ${\bf A}_{1,1}={\bf A}_{2,1}=0$ holds, then $\sigma {\bf A}_{1,1}+
{\bf A}_{2,1}=2(\alpha_5-\alpha_7)$ implies the contradiction $\alpha_5=\alpha_7$.
Therefore ${\bf A}_{1,1}={\bf A}_{2,1}=0$ does not occur,
and moreover if $LR(\chi)\ne0$, then  we take $M\ne0$ and applying the above,  we see that
 there is a non-zero rational vector $(a,b)$
such that vectors  $(c_1,c_3),(c_2,c_4),(-A_{2,1},A_{1,1})$ are proportional  to $(a,b)$.
Taking rational numbers $\kappa_1,\kappa_2$ such that $(c_1,c_3)=\kappa_1(a,b),
(c_2,c_4)=\kappa_2(a,b)$, we have
\begin{equation*}
M=\kappa_1(a{v_1}+b{v_2})+\kappa_2(a\sigma {v_1}+b\sigma {v_2}).
\end{equation*} 
Thus we have proved
\begin{thm}
If $LR(\chi)\ne 0$, then
$(-{\bf A}_{2,1}, {\bf A}_{1,1})$ is proportional  to some non-zero rational vector $(a,b)$,
and $LR(\chi)$ is spanned by $a{v_1}+b{v_2},a\sigma {v_1}+b\sigma {v_2}$.
\end{thm}
\vspace{2mm}

\noindent
{\bf Example}
Let $f(x) = x^8 + 4x^6 + 2x^4 + 28x^2 + 1$ whose Galois group is isomorphic to the dihedral group.
Taking a root $\alpha$, other roots are
\begin{align*}
\alpha_1&=\alpha,
\\
\alpha_2&=-5/24\alpha^7 - 1/24\alpha^6 - 19/24\alpha^5 - 5/24\alpha^4 - 5/24\alpha^3 - 13/24\alpha^2
\\&\hspace{5mm} - 127/24\alpha - 29/24,
\\
\alpha_3&= -5/12\alpha^7 - 19/12\alpha^5 - 5/12\alpha^3 - 139/12\alpha, 
\\
\alpha_4&=-5/24\alpha^7 + 1/24\alpha^6 - 19/24\alpha^5 + 5/24\alpha^4 - 5/24\alpha^3 + 13/24\alpha^2
\\&\hspace{5mm}  - 127/24\alpha  + 29/24,
\\
 \alpha_k&=-\alpha_{k-4}\,\,\,(k=5,6,7,8)
\end{align*}
Automorphisms $\sigma,\mu$ are defined by $\sigma(\alpha_1)=\alpha_2,\mu(\alpha_1)=\alpha_5$.
Then linear relations are spanned by
\begin{align*}
[1]&:\alpha_1+\alpha_2+\alpha_3+\alpha_4+\alpha_5+\alpha_6+\alpha_7+\alpha_8=0\text{ for }\chi_1, 
\\
[2]&:\alpha_1-\alpha_2+\alpha_3-\alpha_4+\alpha_5-\alpha_6+\alpha_7-\alpha_8=0\text{ for }\chi_3, 
\\
[3]&:\alpha_1-\alpha_2+\alpha_3-\alpha_4-\alpha_5+\alpha_6-\alpha_7+\alpha_8=0\text{ for }\chi_4, 
\\
[4]&: 4(v_1+v_2)(\alpha)=\alpha_1-\alpha_2-\alpha_3+\alpha_4+\alpha_5-\alpha_6-\alpha_7+\alpha_8=0\text{ for }\chi_5, 
\\
[5]&: 4\sigma(v_1+v_2)(\alpha)=\alpha_1+\alpha_2-\alpha_3-\alpha_4+\alpha_5+\alpha_6-\alpha_7-\alpha_8=0\text{ for }\chi_5, 
\end{align*}
since we can take $(a,b) = (1,1)$ by ${\bf A}_{1,1}={\bf A}_{2,1}$.
The rank of the matrix $M$ defined by $(\alpha_1,\alpha _2,\dots,\alpha_8)=(1,\alpha,\dots,\alpha^7)M$ is   $3$,
and so there are five relations, which  are intuitively $ \alpha_1-\alpha_2+\alpha_3-\alpha_4=0,\alpha_k=-\alpha_{k-4}$ $(k=5,6,7,8)$.
They correspond  to $([2]+[3])/2$, $([1]+[2]+[4]+[5])/4$,  $([1]-[2]-[4]+[5])/4$,
  $([1]+[2]-[4]-[5])/4$,  $([1]-[2]+[4]-[5])/4$, respectively.  


\end{document}